\documentclass{tsp-short-ukr}

\newtheorem{thm}{Theorem}[section]
\newtheorem{lem}{Lemma}[section]

\theoremstyle{remark}

\newtheorem{dfn}{Definition}[section]
\newtheorem{exm}{Example}[section]
\renewcommand{\proofname}{Proof}

\newcommand{\pr}{\textsf{P}}

\DeclareMathOperator{\M}{\mathsf E}
\usepackage{listings}
\usepackage{color}

\definecolor{dkgreen}{rgb}{0,0.6,0}
\definecolor{gray}{rgb}{0.5,0.5,0.5}
\definecolor{mauve}{rgb}{0.58,0,0.82}
\begin{document}

\title[Maximum of the Brownian Sheet on susubsets]{On the Distribution of Maximum of a Brownian Sheet Restricted to a Lower-Dimensional Set}

\author{N. V.Kruglova}
\address{National Technical University of Ukraine "KPI", Department of Mathematical Analysis and Probability Theory, Pr. Peremohy 37, 02056 Kiev, Ukraine}
\email{natahak@ukr.net}
\author{O. O. Dykhovychnyi}
\address{National Technical University of Ukraine "KPI", Department of Mathematical Analysis and Probability Theory, Pr. Peremohy 37, 02056 Kiev, Ukraine}
\email{a.dyx@ukr.net}

\begin{abstract}
We obtain sufficient conditions of stochastic equivalence of Gaussian random fields with special covariance function. These results generalize Doob's transformation (condition of stochastic equivalence of a Gaussian and a Wiener processes) to the case of random fields (condition of stochastic equivalence of a Gaussian process and a Brownian sheet). We look at the problem of finding the distribution of  supremum of a Brownian sheet on a set with a dimension lower than the dimension of the field. We consider the probability of a Brownian sheet with a certain drift to attain zero level. The obtained results can significantly simplify the problem of finding distributions of functionals of a Brownian sheet by reducing it to the problem of finding distributions on parallelepipeds with dimension lower than the dimension of the field. We consider examples that verify validity of the obtained theorem by modeling corresponding fields and comparing empirical probabilities with theoretical ones.

\end{abstract}

\maketitle

\section{Introduction}
We consider the distribution of an $n$-dimensional Brownian sheet (Chentsov field) $X(s_{1},...,s_{n})$ on a set $S$ of a dimension lower than the dimension of the field. The Brownian sheet was first described by Chentsov \cite{Ch} and Yeh \cite{Y}. In Russian-language literature it was previously known as Chenstov field, whereas in English-language literature a commonly used name was Wiener-Yeh field. Now, it is usually referred to as a Brownian sheet. It is known that random functions such as Brownian sheet occur in modeling of external effects influencing a system at a random moment in time and at a random location. For instance, this is a common scenario in the problem of modeling small transverse vibrations of a string under influence of random external forces or in the problem of heat transfer in a rod with presence of random heating/cooling sources \cite{A},\cite{B} as well as in filtration problems \cite{DKK}.

The problem of finding distribution on a set of a dimension lower than the dimension of the field occurs, for instance, in percolation theory \cite{MMC}.

Although some research on the distribution of maximum of a Brownian sheet on a unit square has been conducted, no exact results have been obtained yet.

For instance, a lot of results regarding the distribution of superemum of a field restricted to curves and polylines have been obtained. For the case of 2-dimensional field on the unit square, Park, Paranjape \cite{PP} \cite{PS} have obtained the distribution of superemum of the field restricted to a polyline with a single vertex. Also, by considering limit of the polyline they have obtained  distribution on the boundary of a square. Generalization of these results for a polyline with $n$ vertices was considered by Klesov and Kruglova \cite{KK} \cite{Kr}. However, the obtained exact distribution of supreremum is given in terms of quite cumbersome integrals direct estimation of which may become problematic. That is why Kruglova and Dykhovychnyi \cite{DK} suggested an empirical approach to finding the distribution by modeling corresponding restriction to polylines. The suggested method is based on Doob's Transformation Theorem \cite{D}.

Also, the problem of finding the following probability is of a special interest:
$$
	\pr\left\{ {\mathop {\sup }\limits_{0\leqslant t \leqslant T} w(t) - f(t) < b} \right\},
$$
where $w(t)$ is a Wiener process, and $f(t)$ is a deterministic function (drift).

For a multidimensional Brownian sheet $X(s_{1},...,s_{n})$, this problem can be formulated as follows:
\begin{equation}
	\label{trivial}
	\pr\left\{\sup_{S}\left(X(s_{1},...,s_{n})-g(s_{1},...,s_{n})\right)<0\right\},
\end{equation}
where $S$ is a set of dimension lower than $n$, and $g(\cdot)$ is a continuous function.

\section{Definitions and preliminaries}

To prove our results the following definitions, lemmas and theorems are needed.
Let $K = [0,1]^n$ denote an $n$-dimensional unit hypercube.
\begin{dfn} A real-valued separable Gaussian stochastic process $\{X\left(s_{1} ,...,s_{n} \right)$, $\overrightarrow{s}=\left(s_{1} ,...,s_{n} \right)\in K \}$ is called a Brownian sheet with $n$ parameters, if $X\left(s_{1} ,...,s_{n} \right)$ is such that
\begin{enumerate}
	\item $X\left(s_{1} ,...,s_{n} \right)=0$, if $s_{1} \cdot ...\cdot s_{n}=0$;
	\item $E\left[X(s_{1} ,...,s_{n} )\right]=0$, for all $(s_{1} ,...,s_{n}) \in K$;
	\item $E\left[X(s_{1} ,...,s_{n} )X(t_{1} ,...,t_{n} )\right]=\prod _{i=1}^{n}\min \{ s_{i} ,t_{i} \}$,  for all $(s_{1} ,...,s_{n} )\in K$ and $(t_{1} ,...,t_{n} )\in K$.
\end{enumerate}
\end{dfn}

\begin{dfn}\label{stekv}
Processes $\xi (t)$ and $\xi^{*} (t)$ are said to be stochastically equivalent if, for all $t \in T$,
$$
	\pr\left\{ {\xi (t) = \xi ^* (t)} \right\} = 1.
$$
\end{dfn}

\begin{thm}\label{D1} (Doob's Transformation Theorem \cite{D}) Let $Y(t)$ be a Gaussian process with $E[Y(t)]=0$, for all $t \ge 0$, and covariance function
\begin{equation}\label{r}
	R(s,t) = u(s)v(t),\qquad s\le t.
\end{equation}
If the function $a(t) = u(t)/v(t)$ is continuous and strictly increasing with the inverse $ a^{-1}(t)$,
then $w(t)$ and $Y(a^{-1}(t))/v(a^{-1}(t))$ are stochastically equivalent processes.
\end{thm}

\section{Main results}

\subsection{Generalization of Doob's Transformation Theorem}
In order to solve (\ref{trivial}), we need the following analogue of Theorem \ref{D1}.
\begin{lem}\label{Dpr1} Let $Y(s_{1} ,..., s_{n})$ be a Gaussian field with $\M\left[Y(s_{1} ,...,s_{n} )\right]=0$, for all $\overrightarrow{s}=(s_{1} ,...,s_{n} )\in [0,\infty )^{n}$, and covariance function
\begin{equation}
	\label{r2}
	\M\left[Y(s_{1} ,...,s_{n} )Y(t_{1} ,...,t_{n} )\right]=\prod _{i=1}^{n}u_{i} (s_{i}\wedge t_{i} )v_{i} (s_{i} \vee t_{i} )
\end{equation}
for all $\overrightarrow{s},\overrightarrow{t}\in [0,\infty )^{n} $ (where $\wedge$ and $\vee$ denote maximum and minimum, respectively).
\\ If the functions $a_{i} (t)=\frac{u_{i} (t)}{v_{i} (t)}$, $i=\overline{1,n}$, are increasing and continuous, then the field
\[
	Z(t_1,\ldots,t_n)=\frac{Y\left(a_{1}^{-1} (t_{1} ),\ldots,a_{n}^{-1} (t_{n} )\right)}{v_{1} \left(a_{1}^{-1} (t_{1} )\right)\cdot \ldots\cdot v_{n} \left(a_{n}^{-1} (t_{n} )\right)}
\]
is stochastically equivalent to $X\left(t_{1} ,\ldots,t_{n} \right)$, where $a_{i}^{-1} (\cdot)$, $i = \overline{1,n}$, denotes the inverse of  $a_{i} (\cdot)$.
\end{lem}
\proofname.

Let us calculate the expectation and the covariance function of $Z(t_1,\ldots,t_n)$. It is clear that $\M\left[Z(t_1,\ldots,t_n)\right]=0$. Then, using (\ref{r2}), we get
$$\M\left[Z(s_1,\ldots,s_n)Z(t_1,\ldots,t_n)\right]=$$
\[\M\left[\frac{Y\left(a_{1}^{-1} (s_{1} ),...,a_{n}^{-1} (s_{n} )\right)}{v_{1} \left(a_{1}^{-1} (s_{1} )\right)\cdot ...\cdot v_{n} \left(a_{n}^{-1} (s_{n} )\right)} \frac{Y\left(a_{1}^{-1} (t_{1} ),...,a_{n}^{-1} (t_{n} )\right)}{v_{1} \left(a_{1}^{-1} (t_{1} )\right)\cdot ...\cdot v_{n} \left(a_{n}^{-1} (t_{n} )\right)} \right]=\]
\[=\frac{u_{1} \left(a_{1}^{-1} (s_{1} )\wedge a_{1}^{-1} (t_{1} )\right)\cdot ...\cdot u_{n} \left(a_{n}^{-1} (s_{n} )\wedge a_{n}^{-1} (t_{n} )\right)}{v_{1} \left(a_{1}^{-1} (s_{1} )\right)\cdot ...\cdot v_{n} \left(a_{n}^{-1} (s_{n} )\right)}\]\[\times\frac{v_{1} \left(a_{1}^{-1} (s_{1} )\vee a_{1}^{-1} (t_{1} )\right)\cdot ...\cdot v_{n} \left(a_{n}^{-1} (s_{n} )\vee a_{n}^{-1} (t_{n} )\right)}{v_{1} \left(a_{1}^{-1} (t_{1} )\right)\cdot ...\cdot v_{n} \left(a_{n}^{-1} (t_{n} )\right)} \]
\[=\frac{u_{1} \left(a_{1}^{-1} (s_{1} )\wedge a_{1}^{-1} (t_{1} )\right)\cdot ...\cdot u_{n} \left(a_{n}^{-1} (s_{n} )\wedge a_{n}^{-1} (t_{n} )\right)}{v_{1} \left(a_{1}^{-1} (s_{1} )\wedge a_{1}^{-1} (t_{1} )\right)\cdot ...\cdot v_{n} \left(a_{n}^{-1} (s_{n} )\wedge a_{n}^{-1} (t_{n} )\right)}\]\[\times\frac{v_{1} \left(a_{1}^{-1} (s_{1} )\vee a_{1}^{-1} (t_{1} )\right)\cdot ...\cdot v_{n} \left(a_{n}^{-1} (s_{n} )\vee a_{n}^{-1} (t_{n} )\right)}{v_{1} \left(a_{1}^{-1} (s_{1} )\vee a_{1}^{-1} (t_{1} )\right)\cdot ...\cdot v_{n} \left(a_{n}^{-1} (s_{n} )\vee a_{n}^{-1} (t_{n} )\right)} \]
\[=a_{1} \left(a_{1}^{-1} (s_{1} )\wedge a_{1}^{-1} (t_{1} )\right)\cdot ...\cdot a_{n} \left(a_{n}^{-1} (s_{n} )\wedge a_{n}^{-1} (t_{n} )\right)\]\[=a_{1} \left(a_{1}^{-1} (s_{1} \wedge t_{1} )\right)\cdot ...\cdot a_{n} \left(a_{n}^{-1} (s_{n} \wedge t_{n} )\right)=(s_{1} \wedge t_{1} ) \cdot ...\cdot  (s_{n} \wedge t_{n} ).\]
The penultimate equality holds because the functions $a_i^{-1}(\cdot)$, $i=\overline{1,n},$ are strictly increasing.
Thus,
$$
	a_{i}^{-1} (s_{i} )\wedge a_{i}^{-1} (t_{i} )=a_{i}^{-1} (s_{i} \wedge t_{i} ), \qquad i=\overline{1,n}.
$$
Using the definition of a Brownian sheet, we get
\[
	\M\left[X(s_{1} ,...,s_{n} )X(t_{1} ,...,t_{n} )\right]=(s_{1} \wedge t_{1} ) \cdot ...\cdot  (s_{n} \wedge t_{n} ),
\]
\[
	\M\left[X(t_{1} ,...,t_{n} )\right]=0.
\]
Due to the fact that the Brownian sheet and $Z(t_1,\cdots,t_n)$ are both Gaussian,
the field $Z(t_1,\cdots,t_n)$ is stochastically equivalent to $X\left(t_{1} ,...,t_{n} \right).$

\begin{lem}\label{Dpr2} Let $Y(s_{1} ,...,s_{n} )$ be a Gaussian field with $\M\left[Y(s_{1} ,...,s_{n} )\right]=0$, for all $\mathop{s}\limits^{\to } =(s_{1} ,...,s_{n} )\in [0,\infty )^{n} $, and covariance function (\ref{r2}). If the functions $a_{i} (t)=\frac{v_{i} (t)}{u_{i} (t)} $, $i=\overline{1,n}$, are decreasing and continuous, then the field
\[
	Z(t_1,\ldots,t_n)=\frac{Y\left(a_{1}^{-1} (t_{1} ),...,a_{n}^{-1} (t_{n} )\right)}{u_{1} \left(a_{1}^{-1} (t_{1} )\right)\cdot ...\cdot u_{n} \left(a_{n}^{-1} (t_{n} )\right)}
\]
is stochastically equivalent to $X\left(t_{1} ,...,t_{n} \right)$, where $a_{i}^{-1} (\cdot)$, $i = \overline{1,n}$, denotes the inverse of  $a_{i} (\cdot)$.
\end{lem}
\proofname .

 Let us calculate the expectation and the covariance function of $Z(t_1,\ldots,t_n)$. Analogously to Lemma \ref{Dpr1}, we get $\M\left[Z(t_1,\ldots,t_n)\right]=0$. Then,
$$\M\left[Z(s_1,\ldots,s_n)Z(t_1,\ldots,t_n)\right]=$$
\[\M\left[\frac{Y\left(a_{1}^{-1} (s_{1} ),...,a_{n}^{-1} (s_{n} )\right)}{u_{1} \left(a_{1}^{-1} (s_{1} )\right)\cdot ...\cdot u_{n} \left(a_{n}^{-1} (s_{n} )\right)} \frac{Y\left(a_{1}^{-1} (t_{1} ),...,a_{n}^{-1} (t_{n} )\right)}{u_{1} \left(a_{1}^{-1} (t_{1} )\right)\cdot ...\cdot u_{n} \left(a_{n}^{-1} (t_{n} )\right)} \right]=\]

\[=\frac{v_{1} \left(a_{1}^{-1} (s_{1} )\vee a_{1}^{-1} (t_{1} )\right)\cdot ...\cdot v_{n} \left(a_{n}^{-1} (s_{n} )\vee a_{n}^{-1} (t_{n} )\right)}{u_{1} \left(a_{1}^{-1} (s_{1} )\vee a_{1}^{-1} (t_{1} )\right)\cdot ...\cdot u_{n} \left(a_{n}^{-1} (s_{n} )\vee a_{n}^{-1} (t_{n} )\right)}\]\[\times\frac{u_{1} \left(a_{1}^{-1} (s_{1} )\wedge a_{1}^{-1} (t_{1} )\right)\cdot ...\cdot u_{n} \left(a_{n}^{-1} (s_{n} )\wedge a_{n}^{-1} (t_{n} )\right)}{u_{1} \left(a_{1}^{-1} (s_{1} )\wedge a_{1}^{-1} (t_{1} )\right)\cdot ...\cdot u_{n} \left(a_{n}^{-1} (s_{n} )\wedge a_{n}^{-1} (t_{n} )\right)} \]

\[=a_{1} \left(a_{1}^{-1} (s_{1} )\vee a_{1}^{-1} (t_{1} )\right)\cdot ...\cdot a_{n} \left(a_{n}^{-1} (s_{n} )\vee a_{n}^{-1} (t_{n} )\right)\]\[=a_{1} \left(a_{1}^{-1} (s_{1} \wedge t_{1} )\right)\cdot ...\cdot a_{n} \left(a_{n}^{-1} (s_{n} \wedge t_{n} )\right)= (s_{1} \wedge t_{1} )\cdot ...\cdot (s_{n} \wedge t_{n} ) .\]
The above holds because the functions $a_i^{-1}(\cdot)$, $i=\overline{1,n},$ are strictly decreasing. Therefore,
$$
	a_{i}^{-1} (s_{i} )\vee a_{i}^{-1} (t_{i} )=a_{i}^{-1} (s_{i} \wedge t_{i} ), \qquad i=\overline{1,n}.
$$
Since
\[
	\M\left[X(t_{1} ,...,t_{n} )X(s_{1} ,...,s_{n} )\right]=(s_{1} \wedge t_{1} ) \cdot ...\cdot  (s_{n} \wedge t_{n} )
\]
and
\[
	\M\left[X(t_{1} ,...,t_{n} )\right]=0,
\]
the field $Z(t_1,\cdots,t_n)$ is stochastically equivalent to $X\left(t_{1} ,...,t_{n} \right).$

\subsection{Main Theorem}
Let $G=[0,y_1]\times\ldots\times [0,y_n]$, where $y_i \in [0,1], i=\overline{1,d}.$

Suppose that $d<n$, and let us consider a $d$-dimensional $S \subset K$ given by
\[\left\{\begin{array}{l} {s_{d+1} =f_{d+1} (s_{j_1}),} \\ {........} \\ {s_{n} =f_{n} (s_{j_d} ).} \end{array}\right., \]
where $j_k\in \{1,\ldots,d\}$, $k=\overline{1,d}$, $0\le s_{i} \le y_i$,  $i=\overline{1,n}$.

Suppose also that there exist decreasing functions $z_{i}(\cdots) ,i=\overline{1,d}$, such that the functions $f_{i+d} (\cdot), i=\overline{1,n-d},$ are decreasing and
\[f_{d+1} (s_{j_1} )\cdot ...\cdot f_{n} (s_{j_d} )=\prod _{i=1}^{d}z_{i} (s_{i} ).\]
 Let $a_{i} (t)=\frac{t}{z_{i} (t)} $ and $x_{i} =a_{i}(y_i)$, $i=\overline{1,d}$, and $D=I_{1} \times ...\times I_{d} $, where $I_{i} =[0,x_{i} ]$.

\begin{thm}\label{rpov}  Let $X(s_{1} ,...,s_{n} )$ be an $n$-parameter Brownian sheet on $G$. Let $g_{S} \left(s_{1} ,...,s_{d} \right)$ denote the restriction of $g\left(s_{1} ,...,s_{n} \right)$ to $S$. If, for all $i=\overline{1,d}$, there exists $a_{i}^{-1} :I_{i} \to [0,y_i]$, the inverse of $a_{i} (\cdot)$,  then
$$\pr\left\{\mathop{\sup }\limits_{S} \left(X\left(s_{1} ,...,s_{n} \right)-g\left(s_{1} ,...,s_{n} \right)\right)< 0\right\}=$$$$=\pr\left\{\mathop{\sup }\limits_{D} \left(X\left(t_{1} ,...,t_{d} \right)-\frac{g_{S} \left(a_{1}^{-1} (t_{1} ),...,a_{d}^{-1} (t_{d} )\right)}{z_{1} \left(a_{1}^{-1} (t_{1} )\right)\cdot ...\cdot z_{d} \left(a_{d}^{-1} (t_{d} )\right)} \right)<0\right\}.$$
\end{thm}
\proofname .

Let $X_{S} (s_{1} ,\ldots,s_{d} )$ denote the restriction of $X(s_{1} ,\ldots,s_{n} )$ to $S$.  The expectation and the covariance function of the field $X_{S} (s_{1} ,\ldots,s_{d} )$ are
$$\M[X_{S} (s_{1} ,\ldots,s_{d} )]=\M[X(s_{1} ,\ldots,s_{d},f_{d+1}(s_{j_1}),\ldots,f_n(s_{j_d}))]=0.$$
Then,
$$\M[X_{S} (s_{1} ,\ldots,s_{d} )X_{S} (t_{1} ,\ldots,t_{d} )]=$$
$$=\M\left[X(s_{1} ,\ldots,s_{d},f_{d+1}(s_{j_1}),\ldots,f_n(s_{j_d}))X(t_{1} ,\ldots,t_{d},f_{d+1}(t_{j_1}),\ldots,f_n(t_{j_d}))\right]$$
$$=(s_1\wedge t_1)\cdot\ldots( s_d\wedge t_d)\cdot( f_{d+1}(s_{j_1})\wedge f_{d+1}(t_{j_1}))\cdot\ldots (f_n(s_{j_d})\wedge f_n(t_{j_d}))$$
$$=(s_1\wedge t_1)\cdot\ldots (s_d\wedge t_d)\cdot (f_{d+1}(s_{j_1}\vee t_{j_1}))\cdot\ldots (f_{n}(s_{j_d}\vee t_{j_d}))$$
$$=\prod_{i=1}^d (s_i\wedge t_i) z_i(s_i\vee t_i).$$
Since $z_{i}(t) ,i=\overline{1,d}$, are decreasing, $\frac{t}{z_i(t)},i=\overline{1,d}$ are increasing. Hence, the assumptions of the theorem allow us to apply Lemma \ref{Dpr1}.
\[\pr\left\{\mathop{\sup }\limits_{S} \left(X\left(s_{1} ,...,s_{n} \right)-g\left(s_{1} ,...,s_{n} \right)\right)< 0\right\}\]\[=\pr\left\{\mathop{\sup }\limits_{[0,y_1]\times\ldots\times[0,y_d] } \left(X_{S} \left(s_{1} ,...,s_{d} \right)-g_{S} \left(s_{1} ,...,s_{d} \right)\right)< 0\right\}\]
\[=\pr\left\{\mathop{\sup }\limits_{D} \left(\frac{X_{S} \left(a_{1}^{-1} (t_{1} ),...,a_{d}^{-1} (t_{d} )\right)}{z_{1} \left(a_{1}^{-1} (t_{1} )\right)\cdot ...\cdot z_{d} \left(a_{d}^{-1} (t_{d} )\right)} -\frac{g_{S} \left(a_{1}^{-1} (t_{1} ),...,a_{d}^{-1} (t_{d} )\right)}{z_{1} \left(a_{1}^{-1} (t_{1} )\right)\cdot ...\cdot z_{d} \left(a_{d}^{-1} (t_{d} )\right)} \right)< 0\right\}\]
\[=\pr\left\{\mathop{\sup }\limits_{D} \left(X\left(t_{1} ,...,t_{d} \right)-\frac{g_{S} \left(a_{1}^{-1} (t_{1} ),...,a_{d}^{-1} (t_{d} )\right)}{z_{1} \left(a_{1}^{-1} (t_{1} )\right)\cdot ...\cdot z_{d} \left(a_{d}^{-1} (t_{d} )\right)} \right)< 0\right\}.\]

Therefore, by using Theorem \ref{rpov}, one can reduce the problem of finding the distribution of maximum of $X\left(t_{1} ,...,t_{d} \right)$ on $S$ to the problem of finding its distribution on the parallelepiped $D$, which may be easier to deal with. For instance, this may be the case if one wants to study the asymptotic behavior of this distribution.

\subsection{Examples}

Consider the two following applications of Theorem \ref{rpov}.
\begin{exm}
Let $X(s_{1} ,s_{2}, s_{3})$ be a 3-dimensional Brownian sheet ($n=3$ and  $d=1$), and
let $S$ be a curve given by
$$
	\left\{\begin{array}{l} {s_{2} =\sqrt{1-s_{1} } ,} \\ {s_{3} =\sqrt{1-s_{1} } .} \end{array}\right.,
$$
Suppose $g(s_{1} ,s_{2} ,s_{3})=s_{1} + s_{2} s_{3}$. Let us show that
$$\pr\left\{\mathop{\sup }\limits_{S} \left(X(s_{1} ,s_{2} ,s_{3} )-g(s_{1},s_{2}, s_{3}) \right)< 0\right\}= 1-e^{-2}.$$
\end{exm}
{\bf Solution.}
From the definition of $S$ we have
\[z_1(s_{1} )=1-s_{1} ,\]

\[a_1(s_{1} )=\frac{s_{1} }{1-s_{1} } .\]
The inverse of $a_1(\cdot)$ is
\[a^{-1} _1(s_{1} )=\frac{s_{1} }{1+s_{1} } .\]
Then,
\[z_1\left(a^{-1}_1 (s_{1} )\right)=\frac{1}{1+s_{1} } \]
and
\[g_{S} (s_{1} )=s_{1} +1-s_{1} =1.\]
Hence,
\[\frac{g_{S} \left(a^{-1}_1 (s_{1} )\right)}{z_1\left(a^{-1}_1 (s_{1} )\right)} =1+s_{1} .\]

\[\pr\left\{\mathop{\sup }\limits_{S} \left(X(s_{1} ,s_{2} ,s_{3} )-s_{1} -s_{2} s_{3} \right)< 0\right\}=\]\[=\pr\left\{\mathop{\sup }\limits_{t\in [0,\infty )} \left(w(t)-1-t\right)< 0\right\}=1-e^{-2} .\]
The above identity is due to the following result \cite{D}:
\begin{equation}\label{wp}
\pr\left\{ {\mathop {\sup }\limits_{0 \leqslant t
< \infty}\left( w(t) - at - b\right)<0} \right\}=1-e^{-2ab}.
\end{equation}

\begin{exm}
Now, suppose $g(s_{1} ,s_{2} ,s_{3} )=s_{1} +s_{2} ^{2} +s_{3} ^{2}$, and let us show that
$$
	\pr\left\{\mathop{\sup }\limits_{S} \left(X(s_{1} ,s_{2} ,s_{3} )-g(s_{1},s_{2}, s_{3}) \right)< 0\right\}= 1-e^{-4}.
$$	
\end{exm}
{\bf Solution.} Since $g(s_{1} ,s_{2} ,s_{3} )=s_{1} +s_{2} ^{2} +s_{3} ^{2} $, we have

\[g_{S} (s_{1} )=s_{1} +1-s_{1} +1-s_{1} =2-s_{1} ,\]

\[\frac{g_{S} \left(a^{-1}_1 (s_{1} )\right)}{z_1\left(a^{-1}_1 (s_{1} )\right)} =2+s_{1} .\]
Finally, using Theorem \ref{rpov}, we get
\[\pr\left\{\mathop{\sup }\limits_{S}\left( X(s_{1} ,s_{2} ,s_{3} )-s_{1} -s_{2} ^{2} -s_{3} ^{2} \right)< 0\right\}=\]\[=\pr\left\{\mathop{\sup }\limits_{t\in [0,\infty )} \left(w(t)-2-t\right)< 0\right\}=1-e^{-4} .\]

It is reasonable to compare obtained probabilities with corresponding empirical estimates resulting from simulation of the $3$-dimensional Brownian sheet. Let us use the algorithm for simulation of a Gaussian processes with special covariance function suggested in \cite{DK}.

The following is R code simulates the Gaussian processes provided in Examples 3.1 and 3.2 and computes both empirical and theoretical distributions of their maximums.
\begin{lstlisting}
> m<-numeric(10^4)
> n<-numeric(10^4)
> for(i in 1:10^4){
+   t<-seq(0,1-1/1000,length.out=1000)
+   vt<-1-t
+   at<-t/(1-t)
+   D<-diff(at)
+   y<-vt*c(0,cumsum(rnorm(999,0,sqrt(D))))  # the field
+   m[i]<-max(y)                             # compute maximum (Example 3.1)
+   n[i]<-max(y-2+t)                         # compute maximum (Example 3.2)
+ }
> m1<-m[m<1]
> length(m1)/10^4
[1] 0.8683                                   # empirical probability (Example 3.1)
> 1-exp(-2)
[1] 0.8646647                                # theoretical probability (Example 3.1)
> n1<-n[n<0]
> length(n1)/10^4
[1] 0.9833                                   # empirical probability (Example 3.2)
> 1-exp(-4)
[1] 0.9816844                                # theoretical probability (Example 3.2)
\end{lstlisting}

It could be seen that the theoretical results are quite close to the empirical ones in both cases (increasing the number of points in the mesh will clearly increase accuracy of the estimate).

Let us consider another example for a $4$-dimensional Brownian sheet.
\begin{exm}
 Let $X(s_1,s_2,s_3,s_4)$ be a $4$-parameter Brownian sheet. For $\lambda>0$, let us consider the following probability:
\begin{equation}\label{ex2}
P_S(X)=P\left\{\mathop{\sup }\limits_{(s_1,s_2,s_3,s_4)\in S} \left(X\left(s_{1},s_2,s_3,s_{4} \right)-\frac{\lambda s_3s_4}{(s_1+s_3)(s_2+s_4)}\right)< 0\right\},
\end{equation}
where $S=\left\{(s_1,s_2,s_3,s_4):0\leq s_i\leq \frac{1}{2},i=\overline{1,4},s_3=1-s_1,s_4=1-s_2\right\}$.

\end{exm}
{\bf Solution.}
Let $X_S$ denote the restriction of the field $X(s_1,s_2,s_3,s_4)$ to $S$. It is clear that $X_S(s_1,s_2)=X(s_1,s_2,1-s_1,1-s_2)$ and

\[R(\bar{s},\bar{t})=E[X_S(\bar{s})X_S(\bar{t})]=\prod_{i=1}^2(s_i\wedge t_i)(1-s_i\vee t_i).\]

Using the notation of Theorem \ref{Dpr1}, we have
\[z_1(s_1)=1-s_1, z_2(s_2)=1-s_2.\]

Therefore,
 \[a_1(s_1)=\frac{s_1}{1-s_1}, a_2(s_2)=\frac{s_2}{1-s_2}\]

 and
 \[ a_1^{-1}(s_1)=\frac{s_1}{1+s_1}, a_2^{-1}(s_2)=\frac{s_2}{1+s_2}, y_1=\frac{1}{2},y_2=\frac{1}{2}, x_1=a_1(y_1)=1,x_2=a_1(y_1)=1.\]

Using Theorem \ref{rpov}, we can rewrite the probability (\ref{ex2}) in the following form:
$$P_S(X)=P\left\{\mathop{\sup }\limits_{(s_1,s_2)\in [0,1/2]^2} \left(X_S\left(s_{1} ,s_2 \right)-\lambda(1-s_1)(1-s_2)\right)< 0\right\}= $$
$$=P\left\{\mathop{\sup }\limits_{(s_1,s_2)\in (0,1]^2}(1+s_1)(1+s_2)X_S\left(\frac{s_{1}}{1+s_1} ,\frac{s_2}{1+s_2} \right)<\lambda\right\}=$$
$$=P\left\{\mathop{\sup }\limits_{(s_1,s_2)\in (0,1]^2}X\left(s_1 ,s_2 \right)<\lambda\right\}.$$
where $X\left(s_1 ,s_2 \right)$ is a $2$-dimensional Brownian sheet. In this way, the initial problem is reduced to the problem of finding the distribution of supremum of the Brownian sheet on a square.

\section*{Conclusion}

In this paper, we obtained a generalization of Doob's Transformation Theorem for multi-dimensional Gaussian random field.

We used this generalization to reduce the problem of finding the distribution of supremum of a $n$-dimensional Brownian sheet restricted to a lower-dimensional set to the problem of finding its distribution on simpler sets (parallelepipeds).

With the help of a special algorithm we performed modeling of a Brownian sheet restrictions to various curves and obtained their empirical distribution which appeared to coincide with the theoretical ones.


\begin{thebibliography}{10}

 \bibitem{A}
 Allen E.J. \textit{Modeling with Ito Stochastic Differential Equations}, Dordrecht: Springer, 2007, 228 p.
doi: 10.1007/978-1-4020-5953-7 .

\bibitem{B}
Bovkun V.A.\textit{ On models that lead to an infinite-dimensional stochastic Cauchy problem}, Teor.
Veroyatnost. i Primenen., 2017, vol. 62, no. 4, pp. 803–804 (in Russian). doi: 10.4213/tvp5151 .

\bibitem{DKK}
Davar Khoshnevisan \textit{Five Lectures on Brownian Sheet}  (2001)  http://www.math.utah.edu/?davar

\bibitem{MMC}
 Меньшиков М. В.,  Молчанов С. А.,  Сидоренко А. Ф. \textit{Теория перколяции и некоторые приложения}, Итоги науки и техн. Сер. Теор. вероятн. Мат. стат. Теор. кибернет., 24, ВИНИТИ, М., 1986, 53–110; J. Soviet Math., 42:4 (1988), 1766–1810]]

\bibitem{D}   Doob J. L. \textit{Heuristic approach to Kolmogorov-Smirnov theorems},
Ann. Math. Statist., \textbf{20}  (1949), 393--403.

\bibitem{Y}
Yeh J. \textit{ Wiener Measure in a Space of Functions of Two Variables} Transactions of the American
Mathematical Society 95 (1960), 433–450.

\bibitem{Ch}
Chentsov N.N. \textit{ Wiener random fields depending on several parameters} Doklady AN SSSR 46
(1956), no. 4, 607–609. (Russian)

\bibitem{M}     Malmquist S.
\textit {On certain confidence contours for distribution functions},
Ann. Math. Statist., \textbf {25}   (1954), 523-533.

\bibitem{PP} S. R. Paranjape and C. Park  \textit{ Distribution of the supremum of the two-parameter Yeh-Wiener process on the boundary}, J. Appl. Probab., \textbf{10} (1973), no~ 4, 875-880.

\bibitem{PS}  Park C. and  Schuurmann F. J.\textit {Evaluations of barrier-crossing probabilities of wiener paths}, J. Appl. Prob., \textbf {13}  (1976),  267-275.
\bibitem{P} Park C.
\textit{ Representations of Gaussian processes by
Wiener processes}, Pacific Journal of Mathematics, \textbf{94} (1981), no~ 2, 407-415.
 \bibitem{Kr}   Kruglova N. V.
\textit{ Distribution of the maximum
of the Chentsov random field},
Theor. Stoch. Process.
 (2008),
no~ 1, 76--81.




\bibitem{KK} Klesov  O.I.,  Kruglova N. V.
\textit {Distribution Of the functionals of  the Chentsov random field.
$R^{3}$}, Naukovi visti NTUU-KPI, \textbf {6} (2007), 145-150 (in ukrainian).
\bibitem{Pro} N. V. Prokhorenko\textit {Stochastic Equivalence of Gaussian Process to the Wiener Process, Brownian Bridge, Ornstein--Uhlenbeck Process}, Naukovi visti NTUU-KPI, \textbf {4}  (2016), 85--93 (in ukrainian).

\bibitem{DK}  Dykhovychnyi O.O., Kruglova N.V.,  Virstiuk O.I.\textit {Methods for identification of probability distribution of random variables from data samples with R statistical computing language}, Mathematics in Modern Technical University, \textbf {1}  (2018), 91-100 (in ukrainian).
\end{thebibliography}
\end{document}